\documentclass[11 pt]{article}

\usepackage{amssymb}
\usepackage{amsmath}

\def\nmark{\mbox{$\rm\bf\kern0.2em\rule{0.06em}{1.45ex}\kern-0.3em
N$}}
\def\dmark{\mbox{$\rm\bf\kern0.2em\rule{0.06em}{1.45ex}\kern-0.3em
D$}}
\def\cmark{\mbox{$\rm\bf\kern0.2em\rule{0.06em}{1.45ex}\kern-0.3em
C$}}
\def\rmark{\mbox{$\rm\bf\kern0.2em\rule{0.06em}{1.45ex}\kern-0.3em
R$}}

\begin{document}
\title{\large \bf
Numerical range of weighted composition operators which contain zero}
 \author{Mahsa Fatehi and Asma Negahdari}

{\maketitle}
\begin{abstract}
In this paper, we  study when zero belongs to   the numerical range of weighted composition operators $C_{\psi,\varphi}$ on the Fock space $\mathcal{F}^{2}$, where $\varphi(z)=az+b$, $a,b \in \mathbb{C}$ and $|a|\leq 1$. In the case that $|a|<1$, we obtain a set contained in the numerical range of $C_{\psi,\varphi}$ and  find the conditions under which   the numerical range of $C_{\psi,\varphi}$ contain zero. Then for $|a|=1$,  we precisely determine the numerical range of $C_{\psi,\varphi}$ and  show that zero lies in its numerical range.
\end{abstract}

\footnote{AMS Subject Classifications. Primary 47B33.\\
{\it Key words and phrases}:  Fock space, Weighted composition operator, Numerical  range.}

\section{Introduction}

The Fock space $\mathcal{F}^{2}$ consists of all entire functions on the complex plane ${\mathbb C}$ which are square integrable with $d\mu(z)=\pi^{-1}e^{-|z|^{2}}dA(z)$ that $dA$ is the Lebesgue measure on ${\mathbb C}$. For $f, g$
in $\mathcal{F}^{2}$, the inner product on the Fock space is given by
$$\langle f,g\rangle=\int_{\mathbb{C}}f(z)\overline{g(z)}d\mu(z).$$
The set $\{e_{m}(z)=z^{m}/\sqrt{m !}:m\geq 0\}$ is an orthonormal basis for $\mathcal{F}^{2}$.
The reproducing kernel at $w$ in $\mathbb{C}$ for $\mathcal{F}^{2}$ is given by $K_{w}(z)=e^{\overline{w}z}$. Let $k_{w}$ denote
the normalized reproducing kernel given by $k_{w} = K_{w}/\|K_{w}\|$,
where $\|K_{w}\|=e^{|w|^{2}/2}$. Fock space is a very important tool for quantum stochastic calculus in the quantum probability. Fore more information about the Fock space, see \cite{z}.\par
Through this paper, for a bounded operator $T$ on  $\mathcal{F}^{2}$, the spectrum of $T$ and the point spectrum of $T$ are denoted by $\sigma(T)$ and $\sigma_{p}(T)$; respectively. For an entire function $\varphi$, the composition operator $C_{\varphi}$ on $\mathcal{F}^{2}$
 is defined by the rule $C_{\varphi}(f)=f \circ \varphi$ for each $f \in \mathcal{F}^{2}$. For an entire function $\psi$, the  weighted composition operator $C_{\psi,\varphi}:\mathcal{F}^{2}\rightarrow \mathcal{F}^{2}$ is given by $C_{\psi,\varphi}h=\psi\cdot(h\circ\varphi)$. There is a vast  literature on composition operators on the other spaces (see  \cite{cm1} and \cite{sh}).  Moreover, recently many authors have worked on the weighted composition operators on the Fock spaces (see \cite{car}, \cite{mf}, \cite{lef},  \cite{zhao1}  \cite{zhao2} and \cite{zhao}).
 Bounded and compact composition operators on the Fock space over $\mathbb {C}^{n}$ were characterized in \cite{car} by Carswell et al. They showed that  $C_{\varphi}$ is bounded on the Fock space if and only if $\varphi(z)=az+b$, where $|a|\leq 1$ and if $|a|=1$, then $b=0$. In \cite{ueki}, Ueki found a necessary and sufficient condition for $C_{\psi,\varphi}$ to be bounded and compact. After that in \cite{lef}, Le gave the easier characterizations for the boundedness and compactness of $C_{\psi,\varphi}$. Moreover, he found normal and isometric weighted composition operators on $\mathcal{F}^{2}$. Unitary  weighted composition operators and their spectrum on the Fock space of $\mathbb {C}^{n}$  were characterized by Zhao in \cite{zhao1}.  Note that  there are some interesting papers  \cite{car},  \cite{zhao1} and \cite{zhao} which were written in another Fock space (see \cite{z}), but their results hold for $\mathcal{F}^{2}$ by the same idea. Then we use them frequently in this paper.\par
 For $T$  a bounded linear operator on a Hilbert space $H$,  the numerical range of $T$ is denoted by    $W(T)$ and is given by $W(T)=\{\langle Tf,f\rangle:\|f\|=1\}$. The set $W(T)$ is convex, its closure contains $\sigma(T)$ and $\sigma_{p}(T) \subseteq W(T)$. There are some  papers that the numerical range of composition operators and weighted composition operators on the Hardy space $H^{2}$ were investigated (see \cite{bs1}, \cite{bs2}, \cite{g2} and \cite{ma}).\par
  In Section 2,   we investigate $W(C_{\psi,\varphi})$, where $\varphi(z)=az+b$ with $0<|a|<1$. In Proposition 2.1, we find a subset contained in  $W(C_{\psi,\varphi})$, where $\psi(\frac{b}{1-a})\neq 0$.  In Theorem 2.2, we show that if $C_{\psi,az+b}$ is compact, where $\psi(\frac{b}{1-a})\neq 0$ and $a$ is not a positive real number, then $W(C_{\psi,\varphi})$ contains zero. Then in Theorem 2.3, for $C_{\psi,\varphi}$ with $\psi(\frac{b}{1-a})=0$, we show that $W(C_{\psi,\varphi})$ contains a closed disk with center at $0$. Moreover, in Remark 2.4, for a constant function $\varphi$, we show that $W(C_{\psi,\varphi})$ contains zero.\par
 In Section 3, for $\varphi(z)=az+b$, with $|a|=1$, we  find the numerical range of $C_{\psi,\varphi}$ and  see that $W(C_{\psi,\varphi})$ contains zero.\\ \par

\section{$\varphi(z)=az+b$, whit $|a|< 1$}

Suppose that $\psi$ is an entire function and  $\varphi(z)=az+b$, where $|a|<1$. If $C_{\psi,\varphi}$ is a bounded operator  on $\mathcal{F}^{2}$, then by \cite[Theorem 1]{zhao2}, $0 \in \sigma(C_{\psi,\varphi})$. Hence, $0 \in \overline{W(C_{\psi,\varphi})}$.  In this section, we study when $0$ belongs to $W(C_{\psi,\varphi})$
and we work on  the numerical range of bounded weighted composition operator $C_{\psi,\varphi}$, where $\varphi(z)=az+b $ with $|a|< 1$. In this section, we assume that $q(z)=e^{\overline{p}(a-1)z}\psi(z+p)$, where $p=\frac{b}{1-a}$ is the fixed point of $\varphi$. In the proof of Proposition 2.1, we will see that $q$ belongs  to  $\mathcal{F}^{2}$ and we assume that $\sum _{j=0}^{\infty}\widehat{q}_{j}\frac{z^{j}}{\sqrt{j!}}$ is the representation series of $q$ in $\mathcal{F}^{2}$.\\ \par

{\bf Proposition 2.1.} {\it Suppose that $\psi$ is an entire function and $\varphi(z)=az+b$, where $0<|a|<1$. Let $C_{\psi,\varphi}$ be bounded on $\mathcal{F}^{2}$. Suppose that  $\psi(p) \neq 0$, where $p=\frac{b}{1-a}$ is the fixed point of $\varphi$.  Let $n$ be a non-negative integer and $m$ be a positive integer. Then $W(C_{\psi,\varphi})$ contains the ellipse with foci at $a^{n}$ and $ a^{n+m}$ and a major axis
 $$\sqrt{|a^{n}-a^{n+m}|^{2}+\frac{|\widehat{q}_{m}a^{n}\sqrt{(m+n)!}|^{2}}{m! n!}}$$
 and a minor axis
 $$\frac{|\widehat{q}_{m}a^{n}|\sqrt{(m+n)!}}{\sqrt{m! n!}}.$$
 }\bigskip

{\bf Proof.}  By \cite[Corollary 1.2]{zhao1},  $C_{k_{p},z-p}$ is unitary and  \cite[Proposition 3.1]{lef} implies that $C_{k_{p},z-p}^{\ast}=C_{k_{-p},z+p}$.
Since
$\varphi(z+p)-p=a(z+p)+b-p=az$
and
$$k_{-p}(z)k_{p}(\varphi(z+p))\psi(z+p)
=e^{\overline{p}(a-1)z}\psi(z+p),$$
we obtain that
\begin{eqnarray*}
C_{k_{p},z-p}^{\ast}C_{\psi,\varphi}C_{k_{p},z-p}&=&C_{k_{-p},z+p}C_{\psi,\varphi}C_{k_{p},z-p}\\
&=&C_{q,az},
 \end{eqnarray*}
where $q=e^{\overline{p}(a-1)z}\psi(z+p)$ (since $C_{q,az}$ is a bounded operator on $\mathcal{F}^{2}$, $q=C_{q,az}(1)$ belongs to $\mathcal{F}^{2}$). It shows that $C_{\psi,\varphi}$ is unitary equivalent to to $C_{q,az}$. Thus, $W(C_{\psi,\varphi})=W(C_{q,az})$ and  so we investigate the numerical range of $C_{q,az}$. Let $M=\mbox{span}\{e_{1},e_{2}\}$, when $e_{1}(z)=\frac{z^{n}}{\sqrt{n!}}$ and $e_{2}(z)=\frac{z^{n+m}}{\sqrt{(n+m)!}}$. We can see that
\begin{eqnarray*}
C_{q,az}(e_{1})(z)&=&(1+\widehat{q}_{1}z+\widehat{q}_{2}\frac{z^{2}}{\sqrt{2!}}+\cdot\cdot\cdot)\frac{a^{n}z^{n}}{\sqrt{n!}}\\
&=&(\frac{a^{n}}{\sqrt{n!}}z^{n}+\widehat{q}_{1}a^{n}\frac{z^{n+1}}{\sqrt{n!}}+\widehat{q}_{2}a^{n}\frac{z^{n+2}}
{\sqrt{2! n!}}+\cdot\cdot\cdot+\widehat{q}_{m}a^{n}\frac{z^{n+m}}{\sqrt{m!n!}}+\cdot\cdot\cdot)
 \end{eqnarray*}
and
\begin{eqnarray*}
C_{q,az}(e_{2})(z)&=&(1+\widehat{q}_{1}z+\widehat{q}_{2}\frac{z^{2}}{\sqrt{2!}}+\cdot\cdot\cdot)\frac{a^{n+m}z^{n+m}}{\sqrt{(n+m)!}}\\
&=&a^{n+m} \frac{z^{n+m}}{\sqrt{(n+m)!}} +\widehat{q}_{1}a^{n+m}\frac{z^{n+m+1}}{\sqrt{(n+m)!}}+\cdot\cdot\cdot.
\end{eqnarray*}
Let $T$ be the compression of $C_{q,az}$ to $M$. Then the matrix representation of $T$ is
$$
\begin{bmatrix}
a^{n} & 0 \\
\widehat{q}_{m}a^{n}\frac{\sqrt{(n+m)!}}{\sqrt{m! n!}} & a^{n+m}
\end{bmatrix}
.$$
By \cite[p.3-4]{gus}, $W(T)$ is an ellipse with foci at $a^{n},a^{n+m}$ and the  minor axis $|a^{n}\widehat{q}_{m}|\frac{\sqrt{(n+m)!}}{\sqrt{m!n!}}$ and the major axis $\sqrt{|a^{n}-a^{n+m}|^{2}+|\widehat{q}_{m} a^{n}\frac{\sqrt{(n+m)!}}{\sqrt{m!n!}}|^{2}}$. Since $W(T)\subseteq W(C_{\psi,\varphi})$, the result follows.\hfill $\Box$ \\ \par


Let $\varphi(z)=az$, where $0< a < 1$. Since $C_{az}$ is normal (see \cite[Theorem 3.3]{lef}), \cite[Theorem 2]{car}, \cite[Proposition 2.6]{zhao} and \cite[Theorem 1.4-4, p.16]{gus} state that $\overline{W(C_{az})}=[0,1]$.  By the Open Mapping Theorem, $0$ is not an eigenvalue for $C_{az}$. Invoking \cite[Theorem 1.5-5, p.20]{gus}, $W(C_{az})=(0,1]$ and so $0 $
dose not belong to $W(C_{az})$.
 In the next theorem,  we prove that $0$ belongs to  the  numerical range of compact weighted composition operator $C_{\psi,az+b}$, where  $a$ is not a positive real number. In the proof of Theorem 2.2, we use the notation  $\mathbb{D}$ which is  the open unit disk in the complex plane $\mathbb{C}$. Moreover,  some ideas of the proof of the next theorem is similar to \cite[Proposition 2.1]{ma}.\\ \par

{\bf Theorem 2.2.} {\it  Suppose that $\psi$ is an entire function and  $\varphi(z)=az+b$, where $0<|a|<1$.  Let $\psi(\frac{b}{1-a}) \neq 0$. Assume that $C_{\psi,\varphi}$ is compact on $\mathcal{F}^{2}$. If $a$ is not a positive real number, then $W(C_{\psi,\varphi})$ contains zero and it is closed.}\bigskip

{\bf Proof.} We know that $W(C_{\psi,\varphi})=\psi(p)W(C_{\frac{\psi}{\psi(p)},\varphi})$, where $p=\frac{b}{1-a}$ is the fixed point of $\varphi$. By \cite[Proposition 2.6]{zhao}, $\sigma(C_{\frac{\psi}{\psi(p)},\varphi})=\{0,1,a,a^{2},...\}$ and by \cite[Theorem 7.1, p. 214]{c1}, $\sigma_{p}(C_{\frac{\psi}{\psi(p)},\varphi})=\{1,a,a^{2},...\}$. Since $\sigma_{p}(C_{\frac{\psi}{\psi(p)},\varphi})\subseteq W(C_{\frac{\psi}{\psi(p)},\varphi})$, the convex hull of some arbitrary  elements of $\sigma_{p}(C_{\frac{\psi}{\psi(p)},\varphi})$ is a subset of $ W(C_{\frac{\psi}{\psi(p)},\varphi})$.
 We claim that there is a set $M$ that $M \subseteq W(C_{\frac{\psi}{\psi(p)},\varphi})$ and $0\in M$. We break the problem  into three cases.\\
(a) Assume  $a=|a|e^{i\theta}$ and $e^{i\theta}$ is not a root of $1$. Then $\{e^{i n \theta}: n\geq 0\}$ is dense in $\partial \mathbb{D} $. We can find $n_{1},n_{2},n_{3},n_{4}$ such that $|a|^{n_{1}}e^{i  n_{1}\theta}, |a|^{n_{2}}e^{i n_{2}\theta}, |a|^{n_{3}}e^{i n_{3}\theta}, |a|^{n_{4}}e^{i n_{4}\theta}$ lie in the quadrants I, II, III, IV; respectively. It is not hard to see that $0$ is contained in the interior of the polygonal region $P$ whose vertices are $|a|^{n_{1}}e^{i  n_{1}\theta}, |a|^{n_{2}}e^{i n_{2}\theta}, |a|^{n_{3}}e^{i n_{3}\theta}, |a|^{n_{4}}e^{i n_{4}\theta}$. Let $M$ be the union of $P$ and its interior region. \\
(b) Assume that $a=|a|e^{i \theta}$ that $e^{i \theta}$ is a primitive root of $1$ of order $n >2$. Let $P$ be the polygonal region whose vertices are $1, |a|e^{i \theta}, |a|^{2}e^{2i \theta}, ..., |a|^{n-1}e^{(n-1)i  \theta}$ (note that $e^{i\theta},e^{2i\theta},...,e^{(n-1)i\theta}$ are the $n$th root of $1$). Since $n> 2$, the argument of $a$ is not $0$ or $\pi$ and so there are at least three
vertices which are non-colinear points. It is not hard to see that   none of sides of $P$ contains zero and so $0 $ belongs to the interior of the polygonal region $ P$.   Again let $M$ be the union of $P$ and its interior region. \\
(c) Assume that $a=-|a|$. As we know, the convex hull of $\sigma_{p}(C_{\frac{\psi}{\psi(p)},\varphi})$ is a subset of $W(C_{\frac{\psi}{\psi(p)},\varphi})$. Then $[a,1]\subseteq W(C_{\frac{\psi}{\psi(p)},\varphi})$. Let $M$ be the closed line segment with end points $a$  and $1$.\\
Since in these three cases, $M\subseteq W(C_{\frac{\psi}{\psi(p)},\varphi})$ and $0 \in M$,   $0 \in W(C_{\psi,\varphi})$.  Moreover, invoking \cite[Theorem 1]{de}, $W(C_{\psi,\varphi})$ is closed.  \hfill $\Box$ \\ \par

Note that if $\varphi$ and $\psi$ satisfy the hypotheses of Theorem 2.2 and the argument of $a$ is not $0$ or $\pi$, then by the proof of Theorem 2.2, $0$ lies in the interior of $W(C_{\psi,\varphi})$. In the next theorem, we show that $0$ belongs to the interior of $W(C_{\psi,az+b})$, where $0< |a| <1$ and $\psi(\frac{b}{1-a})=0$.\\ \par

{\bf Theorem 2.3.} {\it  Suppose that $\psi$ is an entire function and $\varphi(z)=az+b$, where $0<|a|<1$. Let $n$ be a non-negative integer and $m$ be a positive integer. Assume that $C_{\psi,\varphi}$ is bounded on $\mathcal{F}^{2}$. Suppose that  $\psi(p) = 0$, where $p=\frac{b}{1-a}$ is the fixed point of $\varphi$.  Then $W(C_{\psi,\varphi})$ contains a closed disk with center at $0$ and radius $|\widehat{q}_{m}a^{n}\frac{\sqrt{(n+m)!}}{m! n!}|/2$.}\bigskip

{\bf Proof.} As we saw in the proof of Proposition 2.1, $W(C_{\psi,\varphi})=W(C_{q,az})$, so we investigate the numerical range of $C_{q,az}$. We assume that $M=\mbox {span}\{e_{1},e_{2}\}$, where $e_{1}(z)=\frac{z^{n}}{\sqrt{n!}}$ and $e_{2}(z)=\frac{z^{n+m}}{\sqrt{(n+m)!}}$.
We have
\begin{eqnarray*}
C_{q,az}(e_{1})&=&(\widehat{q}_{1}z+\widehat{q}_{2}\frac{z^{2}}{\sqrt{2!}}+\cdot\cdot\cdot)a^{n}\frac{z^{n}}{\sqrt{n!}}\\
&=& \widehat{q}_{1}\frac{a^{n}}{\sqrt{n!}}z^{n+1}+\widehat{q}_{2}a^{n}\frac{z^{n+2}}{\sqrt{2! n!}}+\cdot\cdot\cdot+\widehat{q}_{m}a^{n}\frac{z^{n+m}}{\sqrt{n! m!}}+\cdot\cdot\cdot
\end{eqnarray*}
and
\begin{eqnarray*}
C_{q,az}(e_{2})&=&(\widehat{q}_{1}z+\widehat{q}_{2}\frac{z^{2}}{\sqrt{2!}}+\cdot\cdot\cdot)a^{n+m}\frac{z^{n+m}}{\sqrt{(n+m)!}}\\
&=&\widehat{q}_{1}a^{n+m}\frac{z^{n+m+1}}{\sqrt{(n+m)!}}+\cdot\cdot\cdot.
\end{eqnarray*}
Let $T$ be the compression of $C_{\psi,\varphi}$ to $M$. Then the matrix representation of $T$ is \\ \par
$$
\begin{bmatrix}
0 & 0 \\
\widehat{q}_{m}a^{n}\frac{\sqrt{(n+m)!}}{\sqrt{n! m!}} & 0
\end{bmatrix}
.$$
By \cite[Example 1, p. 1]{gus},
$$W(T)=\left\{z: |z|\leq \frac{|\widehat{q}_{m}a^{n}\sqrt{(n+m)!}|}{2\sqrt{n! m!}}\right\}.$$
Since $W(T) \subseteq W(C_{\psi,\varphi}) $, $W(C_{\psi,\varphi})$ contains a closed disk with center at $0$ and radius $\left| \frac{\widehat{q}_{m}a^{n}\sqrt{(n+m)!}}{2\sqrt{n! m!}}\right|$.\hfill $\Box$ \\ \par

 Remark 2.4.  Suppose that for some complex number $b$, $\varphi\equiv b$ and $\psi$  is a non-zero entire function. Then $C_{\psi,\varphi}f=f(b)\psi =\langle f,\|\psi\|K_{b}\rangle \frac{\psi}{\|\psi\|}$. By \cite[Proposition 2.5]{bs2}, we can find $W(C_{\psi,\varphi})$ as follows.\\
 (a) If $K_{b}=\frac{c}{\|\psi\|^{2}}\psi$ for some non-zero complex number $c $, then $W(C_{\psi,\varphi})$ is the closed line segment from $0$ to $\overline{c}$.\\
 (b) If $K_{b}\perp \psi$, then $W(C_{\psi,\varphi})$ is the closed disk  centered at the origin with radius $\frac{\|\psi\|}{2}e^{\frac{|b|^{2}}{2}}$.\\
 (c) Otherwise $W(C_{\psi,\varphi})$ is a closed ellipse  with foci at $0$ and $ \psi(b)$.\\
 Then we can see that in the case that $\varphi$ is constant, $W(C_{\psi,\varphi})$ contains zero. \\ \par
In the first part of the following example, we give a compact weighted composition operator $C_{\psi,az+b}$, where $a$ is a positive real number and $0 \in W(C_{\psi, az+b})$ (see Theorem 2.2). Also in the second part, we give an example which satisfy the conditions of Theorem 2.3.\\ \par

Example 2.5. (a) Suppose that $\varphi(z)=\frac{1}{2}z-\frac{1}{2}$ and $\psi(z)=e^{z}$. By \cite[Corollary 2.4]{zhao}, $C_{\psi,\varphi}$ is compact. It is easy to see that $1$ is the fixed point of $\varphi$ and $q(z)=e K_{1/2}(z)$. The representation series of $q$ in $\mathcal{F}^{2}$ is $$\sum_{j=0}^{\infty}\frac{e}{2^{j}\sqrt{j!}}\frac{z^{j}}{\sqrt{j!}}.$$
Let $n=m=1$. By Proposition 2.1, $W(C_{\psi,\varphi})$ contains the  ellipse with foci at $1/2$ and $1/4$ and the major axis $\sqrt{\frac{1}{16}+\frac{e^{2}}{8}}$. It states that $0$ belongs to $W(C_{\psi,\varphi})$.\par
(b) Let $\varphi(z)=\frac{1}{2}z+\frac{1}{2}$ and $\psi(z)=K_{1}(z)-e^{-1}$. Note that $C_{\psi,\varphi}=C_{K_{1},\varphi}-e^{-1}C_{\varphi}$ and so by \cite[Theorem 2]{car} and \cite[Proposition 2.2]{zhao}, $C_{\psi,\varphi}$ is bounded. We have $\psi(-1)=0$ and $q(z)=e^{z/2}(e^{z-1}-e^{-1})=e^{-1}(e^{\frac{3}{2}}z-e^{\frac{1}{2}z}).$ It is not hard to see that the  representation series of $q$ in $\mathcal{F}^{2}$
is
$$\sum_{j=0}^{\infty}\frac{3^{j}-1}{e 2^{j}\sqrt{j!}}\frac{z^{j}}{\sqrt{j!}}.$$
 Let $n=0$ and $m=1$. By Theorem 2.3, $W(C_{\psi,\varphi})$ contains the closed disk with center at $0$ and radius $\frac{1}{2e}$.\\ \par

\section{$\varphi(z)=az+b$, whit $|a|= 1$}
In this section, we completely find the numerical range of $C_{\psi,az+b}$, where $|a|=1$.
Let $S$ be a subset of complex plane $\mathbb{C}$. For $a \in \mathbb{C} $, we define $aS=\{as: s \in S\}$; we use this definition in the next theorem.  \\ \par

{\bf Theorem 3.1.} {\it Suppose that $\varphi(z)=az+b$, where $|a|=1$. Let  $C_{\psi,\varphi}$ be a bounded weighted composition operator on $\mathcal{F}^{2}$. Then \\
(a) If $a \neq 1$ and $a$ is a primitive root of $1$ of order $n$, then $W(C_{\psi,\varphi})=\psi(0)e^{\frac{a|b|^{2}}{a-1}}P$, where $P$ is the union of the  polygon with $n$ sides and vertices at $1, a, ..., a^{n-1}$ and its interior region.\\
(b) If  $a$ is not a  root of $1$, then $W(C_{\psi,\varphi})=\psi(0)e^{\frac{a|b|^{2}}{a-1}}\mathbb{D} \cup \{\psi(0)e^{\frac{a|b|^{2}}{a-1}}a^{m}:m\geq 0\}$.\\
(c) If $a=1$,  then $W(C_{\psi,\varphi})=\psi(0)e^{\frac{|b|^{2}}{2}}\mathbb{D}$}.\bigskip

{\bf Proof.} Suppose that $a \neq 1$. By \cite[Proposition 2.1]{lef}, $\psi(z)=\psi(0)e^{-a\overline{b}z}=\psi(0)K_{-\overline{a}b}(z)$. Let $u=\frac{-\overline{a}b}{\overline{a}-1}$. We have

\begin{eqnarray}
C_{k_{u},z-u}C_{\psi,\varphi}C_{k_{-u},z+u}&=&\frac{1}{\|K_{u}\|^{2}}C_{e^{\overline{u}z},z-u}C_{\psi,\varphi}C_{e^{-\overline{u}z},z+u}\nonumber\\
&=&\frac{1}{\|K_{u}\|^{2}}e^{\overline{u}z}\cdot \psi(z-u)\cdot e^{(-\overline{u}(az+b))\circ(z-u)}C_{(z+u)\circ (az+b)\circ (z-u)}\nonumber\\
&=&\frac{1}{\|K_{u}\|^{2}}e^{\overline{u}z}\cdot \psi(z-u)\cdot e^{-\overline{u}(az-au+b)}C_{az+u(1-a)+b}\nonumber\\
&=& C_{\widetilde{\psi},\widetilde{\varphi}},
 \end{eqnarray}
where
$$\widetilde{\varphi}(z)=az+\frac{-\overline{a}b}{\overline{a}-1}(1-a)+b=az$$
 and
$$\widetilde{\psi}(z)=e^{-|u|^{2}}e^{\overline{u}z}\cdot \psi(z-u)\cdot e^{-\overline{u}(az-au+b)}=\psi(0)e^{\frac{a|b|^{2}}{a-1}}.$$
Then $C_{\psi,\varphi}$ is unitary equivalent to $\psi(0)e^{\frac{a|b|^{2}}{a-1}}C_{az}$ (see \cite[Corollary 1.2]{zhao1}). We try to find $W(C_{az})$.
 We prove that $\sigma_{p}(C_{az})=\{1, a,a^{2}, ...\}$. It is easy to see that $C_{az}(z^{j})=a^{j}z^{j}$ for each non-negative integer $j$. Then $\{1,a,...\}\subseteq \sigma_{p}(C_{az})$. Since by \cite[Lemma 2]{car}, $C_{az}^{\ast}=C_{\overline{a}z}$,  $C_{az}$ is an isometry. We infer from \cite[Exercise 7, p.213]{c1} that $\sigma_{p}(C_{az}) \subseteq \partial \mathbb{D}$ (note that $C_{az}  $ is invertible). Assume that there is $\lambda \in \sigma_{p}(C_{az})$ such that $|\lambda|=1$ and $\lambda  $ does not belong to $\{a^{m}: m\geq 0\}$. Thus, there exists a non-zero function $f \in \mathcal{F}^{2}$ that
\begin{eqnarray}
C_{az}(f)&=&\lambda f.
\end{eqnarray}
 It shows that $f(0)=\lambda f(0)$. Hence $f(0)=0$. Assume that for each $j < k$, $f^{(j)}(0)=0$. We prove that $f^{(k)}(0)=0$. Taking $k$th derivatives on the both sides of  Equation (2) yields $a^{k}f^{(k)}(0)=\lambda f^{(k)}(0)$. Then $f^{(k)}(0)=0$. Thus, $f\equiv 0$ which is a contradiction. It states that $\sigma_{p}(C_{az})=\{1, a,a^{2},...\}$. Moreover, by \cite[Corollary 1.4]{zhao1}, $\sigma(C_{az})=\overline{\{a^{m}\}}_{m=0}^{\infty}$.\\
 (a) Suppose that $a \neq 1$ and $a$ is a primitive root of $1$ of order $n$. If $n=2$, then $\sigma(C_{az})=\{-1,1\}$.  Invoking \cite[Theorem 1.4-4, p.16]{gus}, $\overline{W(C_{az})}$ is the convex hull of $\sigma(C_{az})$ which is equal to $[-1,1]$. Since $-1,1 \in \sigma_{p}(C_{az})$, we conclude that $W(C_{az})=[-1,1]$. Now assume that $n> 2$. Let $P$ be the convex hull of $\{1,a,...,a^{n-1}\}$ which is the union of  polygon with $n$ sides and vertices at $1, a, ...,a^{n-1}$ and its interior region. We can see that $P$ is a subset of $W(C_{az})$ (note that $\sigma_{p}(C_{\psi,\varphi})= \{1,a,...,a^{n-1}\}$ and $\sigma_{p}(C_{az})\subseteq W(C_{az})$). We show that $W(C_{az})=P$. By \cite[Theorem 1.4-4, p.16]{gus}, $\overline{W(C_{az})}=P$. Since  all vertices of $P$ belong to $\sigma_{p}(C_{az})$,  $W(C_{az})=P$ (see also \cite[Corollary 1.5-7, p.20]{gus}). It shows that $W(C_{\psi,\varphi})=\psi(0)e^{\frac{a|b|^{2}}{a-1}}P$.\\
 (b) Assume that $a$ is not a root of $1$. Since $\sigma_{p}(C_{az})=\{a^{m}\}_{m=0}^{\infty}$ and $\{a^{m}:m\geq 0\}$ is a dense subset of the unit circle, we get $\{a^{m}:m\geq 0\}\cup \mathbb{D}$ is a subset of $W(C_{az})$. Moreover, since $C_{az}$ is an isometry, $\|C_{az}\|=1$ and so $W(C_{az})\subseteq \overline{\mathbb{D}}$. Now we show that for each $\lambda \in \partial \mathbb{D}$ that $\lambda \notin \{a^{m}:m\geq 0\}$, $\lambda \notin W(C_{az})$. Suppose that there is $\lambda \in \partial \mathbb{D}$ which does  not belong to $\{a^{m}: m\geq 0\}$ and $\lambda \in W(C_{az})$. By \cite[Theorem 1.3-3, p.10]{gus}, $\lambda \in \sigma_{p}(C_{az})$ which is a contradiction. Then $W(C_{az})=\{a^{m}:m \geq 0\} \cup \mathbb{D}$ and so $W(C_{\psi,\varphi})=\psi(0)e^{\frac{a|b|^{2}}{a-1}}\mathbb{D} \cup \{\psi(0)e^{\frac{a|b|^{2}}{a-1}}a^{m}:m\geq 0\}$.\\
(c) Assume that $a=1$. By \cite[Proposition 2.1]{lef}, $\psi(z)=\psi(0)K_{-b}(z)$. We know that $C_{\frac{K_{-b}}{\|K_{-b}\|},z+b}$ is unitary (see \cite[Corollary 1.2]{zhao1}). Then $W(C_{\psi,z+b})=\psi(0)\|K_{-b}\|W(C_{\frac{K_{-b}}{\|K_{-b}\|},z+b})$. We try to find $W(C_{\frac{K_{-b}}{\|K_{-b}\|},z+b})$. Since $\sigma(C_{\frac{K_{-b}}{\|K_{-b}\|},z+b})=\partial  \mathbb{D}$ (see \cite[Corollary 1.4]{zhao1}), \cite[Theorem 1.4-4, p. 16]{gus} implies that $\overline{W(C_{\frac{K_{-b}}{\|K_{-b}\|},z+b})}=\overline{\mathbb{D}}$. Since $W(C_{\frac{K_{-b}}{\|K_{-b}\|},z+b})$ is convex, it is not hard to see that $\mathbb{D} \subseteq W(C_{\frac{K_{-b}}{\|K_{-b}\|},z+b})$. Because $C_{\frac{K_{-b}}{\|K_{-b}\|},z+b}$ is unitary, $\|C_{\frac{K_{-b}}{\|K_{-b}\|},z+b}\|=1$. Hence by \cite[Theorem 1.3-3, p.10]{gus}, if $\lambda \in \partial \mathbb{D}$ is an element of $W(C_{\frac{K_{-b}}{\|K_{-b}\|},z+b})$, then $\lambda$ must belong to $\sigma_{p}(C_{\frac{K_{-b}}{\|K_{-b}\|},z+b})$. We claim that $\sigma_{p}(C_{\frac{K_{-b}}{\|K_{-b}\|},z+b})=\emptyset$. Assume that $\lambda \in \sigma_{p}(C_{\frac{K_{-b}}{\|K_{-b}\|},z+b})$. Let $\mu$ be an arbitrary unimodular number. It is not hard to see that there exists $u \in \mathbb{C} $ such that $e^{2\mbox{Im}(u\overline{b})i}=\mu$. By Equation (1), we get
$$C_{k_{u}, z-u}C_{\frac{K_{-b}}{\|K_{-b}\|},z+b}C_{k_{-u}, z+u}=C_{\widetilde{\psi},z+b},$$
where $\widetilde{\psi}(z)=\mu \psi(z)$. Then $C_{\frac{K_{-b}}{\|K_{-b}\|},z+b}$ is unitary equivalent to $\mu C_{\psi,\varphi}$. It states that $\lambda \overline{\mu}$ is an eigenvalue of $C_{\psi,\varphi}$. Thus, $\sigma_{p}(C_{\psi,\varphi})=\partial \mathbb{D}$. Since $C_{\psi,\varphi}$ is normal (see \cite[Theorem 3.3]{lef}) and the Fock space is sparable, by \cite[Proposition 5.7, p.47]{c1}, $C_{\psi,\varphi}$ cannot have an uncountable collection of eigenvalues which is a contradiction. Therefore, $\sigma_{p}(C_{\frac{K_{-b}}{\|K_{-b}\|},z+b})=\emptyset$ and it shows that $W(C_{\frac{K_{-b}}{\|K_{-b}\|},z+b})=\mathbb{D}$. Thus, $W(C_{\psi,\varphi})=\psi(0)e^{\frac{|b|^{2}}{2}}\mathbb{D}$.\hfill $\Box$ \\ \par
In Theorem 3.1, we saw that $0$ lies in the interior of the numerical range of $C_{\psi,az+b}$, where $|a|=1$. In the next example, we compute the numerical range of some weighted composition operators using Theorem 3.1.\\ \par

Example 3.2. (a) Let $\varphi(z)=iz+3$ and $\psi(z)=K_{3i}(z)$. By Theorem 3.1(a), $W(C_{\psi,\varphi})=e^{\frac{9i}{i-1}}P$, where $P$ is the union of the  polygon with $n$ sides and vertices at $1, i, -1,-i$ and its interior region.\par
(b) Let $\varphi(z)=e^{\sqrt{3}i}z+2$ and $\psi(z)=K_{-2e^{-\sqrt{3}i}}(z)$. Theorem 3.1(b) implies that $W(C_{\psi,\varphi})=t\mathbb{D} \cup \{t (e^{\sqrt{3}i})^{m}:m\geq 0\}$, where
$$t=e^{\frac{4e^{\sqrt{3}i}}{e^{\sqrt{3}i}-1}}.$$ \par
(c) Let $\varphi(z)=z+2$ and $\psi(z)=K_{-2}(z)$. We infer from Theorem 3.1(c) that $W(C_{\psi,\varphi})=e^{2}\mathbb{D}$.\\ \par

\bigskip
{M. Fatehi, Department of Mathematics, Shiraz Branch, Islamic Azad
University, Shiraz, Iran. \par E-mail: fatehimahsa@yahoo.com \par
A. Negahdari, Department of Mathematics, Shiraz Branch, Islamic Azad
University, Shiraz, Iran. \par E-mail:  asma.negahdari9219@gmail.com \par
}

\end{document}